\documentclass[11pt]{article}
\usepackage[font=small,labelfont=bf]{caption}
\usepackage{latexsym,amsfonts,amsthm,amsmath,amssymb}
\usepackage{hyperref}
\usepackage{xcolor}
\usepackage{verbatim}
\usepackage{enumerate}
\usepackage{tabularx,graphicx}

\begin{document}

\begin{center}
{\LARGE Euler's First Proof of Stirling's Formula} \\

\bigskip\bigskip

{Alexander Aycock, Johannes-Gutenberg University Mainz \\ Staudinger Weg 9, 55128 Mainz \\ \url{aaycock@students.uni-mainz.de}}

\bigskip

%{Maria Agnesi, University of Bologna \\ 86 Via del Pallone, Bologna, Italy 40100 \\ {\tt magnesi@bologna.edu}}

\end{center}

\medskip

\begin{abstract}\noindent % Abstract not required for review articles.
We present a proof given by Euler in his paper {\it ``De serierum determinatione seu nova methodus inveniendi terminos generales serierum"} \cite{E189} (E189:``On the determination of series or a new method of finding the general terms of series") for Stirling's formula. Euler's proof uses his theory of difference equations with constant coefficients. This theory outgrew from his earlier considerations on inhomogeneous differential equations with constant coefficients of finite order that he tried to extend to the case of infinite order.
\end{abstract}

\medskip

%%%%%%%%%%%%%%%%%%%%%%%%%%%%%%%%%%%%%%%%%%%%%%%%%%%%%%%%%%%%%%%%%%%%%%%%%%%%%%%%
\section{Introduction}
\label{sec: Introduction}

Stirling's formula

\begin{equation}
\label{eq: Stirling}
    n! \sim \sqrt{2\pi n}\left(\dfrac{n}{e}\right)^n \quad \text{for} \quad n \rightarrow \infty.
\end{equation}
was first proven by Stirling. It can be proven by application of the Euler-Maclaurin summation formula or the saddle point approximation. But in his paper {\it``De serierum determinatione seu nova methodus inveniendi terminos generales serierum"} \cite{E189} (E189:``On the determination of series or a new method of finding the general terms of series") Euler gave another proof based on his theory on inhomogeneous linear difference equations with constant coefficients. His theory will be described in section \ref{sec: Euler's Theory of inhomogeneous difference equations  with constant Coefficients}. Finally, we will present and discuss Euler's proof in section \ref{sec: Application to the Factorial}.

%%%%%%%%%%%%%%%%%%%%%%%%%%%%%%%%%%%%%%%%%%%%%%%%%%%%%%%%%%%%%%%%%%%%%%%%%%%%%%%%
\section{Euler's Theory of Inhomogeneous Difference Equations  ith Constant Coefficients}
\label{sec: Euler's Theory of inhomogeneous difference equations  with constant Coefficients}

In this section we will discuss  Euler's application  of his theory of inhomogeneous  difference equations  with constant coefficients to the derivation of Stirling's formula (\ref{eq: Stirling}). Euler reduced them to a differential equation of infinite order. Having treated the finite order case in {\it ``Methodus aequationes differentiales altiorum graduum integrandi ulterius promota"} \cite{E188} (E188:``The method to integrate differential equations of higher degrees expanded further") before, in \cite{E189} he then tried to transfer the results from the before-mentioned paper to the case of infinite order. Unfortunately, this is not possible in the way Euler intended and hence lead Euler to a wrong result when he applied his theory to the case of the logarithm of the factorial. We will explain this in more detail in section \ref{sec: Application to the Factorial}. But we will briefly state what we need to discuss Euler's solution of inhomogeneous linear differential equations of finite (see section \ref{subsec: Inhomogeneous linear Differential Equations of finite Order}) and infinite order (see section \ref{subsec: Reduction of the Difference Equation to a Differential Equation}) first.

\subsection{Inhomogeneous Linear Differential Equations of Finite Order}
\label{subsec: Inhomogeneous linear Differential Equations of finite Order}

In his paper \cite{E188}, Euler considered equations of the form:

\begin{equation}
    \label{eq: Inhomogeneous Finite}
    \left(a_0+a_1 \dfrac{d}{dx}+a_2 \dfrac{d^2}{dx^2} +\cdots+ a_n \dfrac{d^n}{dx^n} \right)f(x)= g(x),
\end{equation}
with complex coefficients $a_1, a_2, \cdots, a_n$. Euler did not state any conditions on the function $g(x)$\footnote{The conditions on $g(x)$ can be inferred from Euler's solution. But since we will not need this in this paper, we will not elaborate on this subject.}. In \S 22, Euler described the following procedure: First, find the zeros with their multiplicity of the expression:

\begin{equation*}
    P(z)=a_0+a_1 z+a_2 z^2 +\cdots+ a_n z^n.
\end{equation*}
Assume $z=k$ is a solution of $P(z)=0$. Then, if $k$ is a simple zero\footnote{In this note, we will only need the case of simple zeros and hence will only state the corresponding formula. In \cite{E188}, Euler stated all cases from order 1 to 4 explicitly.} of $P(z)$, a solution of (\ref{eq: Inhomogeneous Finite}) is given by:

\begin{equation}
    \label{eq: Inhomogeneous Finite Solution}
    f(x)= \dfrac{e^{kx}}{P'(k)} \int e^{-kx}g(x)dx.
\end{equation}
Note that the indefinite integral introduces a constant of integration.

\subsection{Reduction of the Difference Equation to a Differential Equation}
\label{subsec: Reduction of the Difference Equation to a Differential Equation}

\subsubsection{General Idea}

As we mentioned in section \ref{sec: Introduction}, Euler's paper \cite{E189} is a paper actually devoted to inhomogeneous difference equations with constant coefficients, i.e., equations of the form:

\begin{equation}
    \label{eq: Difference equation}
    a_0 f(x)+a_1f(x+1)+\cdots+ a_n f(x+n) = g(x),
\end{equation}
with complex coefficients $a_0, a_1, \cdots, a_n$. Euler's idea to solve (\ref{eq: Difference equation}) is as follows: First, rewrite $f(x+1), f(x+2), \cdots, f(x+n)$ in terms of $f(x)$ and its derivatives by applying Taylor's theorem. Next, substitute the corresponding term in equation (\ref{eq: Difference equation}). After some rearrangement, one arrives at an inhomogeneous differential equation  of infinite order with constant coefficients, i.e., an equation of the form:

\begin{equation}
\label{eq: Differential infinite}
    \left(A_0 +A_1 \dfrac{d}{dx} +A_2 \dfrac{d^2}{dx^2} +\cdots + A_n\dfrac{d^n}{dx^n}+ \cdots \right) f(x) = g(x),
\end{equation}
where $A_0, A_1, A_2, \cdots$ are complex coefficients.

Having transformed the initial equation (\ref{eq: Difference equation}) into this form, Euler argued that the same procedure outlined in section (\ref{subsec: Reduction of the Difference Equation to a Differential Equation}) also applies here. More precisely,  one  has to find all zeros of the expression:

\begin{equation}
\label{eq:  Infinite Polynomial}
    A_0 +A_1 z +A_2 z^2 +\cdots + A_n z^n +\cdots
\end{equation}
and has to construct the solution to (\ref{eq: Differential infinite}) from those zeros. 
In his paper \cite{E189} Euler considered various examples; but in this note we are interested in his solution of the simple difference equation.

\subsubsection{Example: The Simple Difference Equation}

For the sake of explanation and since we will be need the result in section (\ref{sec: Application to the Factorial}), let us consider the simple difference equation, i.e., the equation

\begin{equation}
    \label{eq: Simple Difference}
    f(x+1)-f(x)=g(x)
\end{equation}
and let us describe Euler's solution. First, Euler\footnote{In his paper \cite{E189} \S 55, Euler considered the equation $y(x)-y(x-1)=X(x)$ instead of equation (\ref{eq: Simple Difference}). But does not change the final result substantially, of course.} expanded $f(x+1)$ by using Taylor's theorem:

\begin{equation*}
    f(x+1) = f(x)+\dfrac{d}{dx}f(x)+\dfrac{1}{2!}\dfrac{d^2}{dx^2}f(x)+\dfrac{1}{3!}\dfrac{d^3}{dx^3}f(x) \cdots .
\end{equation*}
Substituting this into equation (\ref{eq: Simple Difference}), Euler arrived at the equation

\begin{equation*}
    \left(\dfrac{d}{dx}+\dfrac{1}{2!}\dfrac{d^2}{dx^2}+ \dfrac{1}{3!}\dfrac{d^3}{dx^3}+\cdots\right)f(x)=g(x).
\end{equation*}
Thus, according to his theory, Euler needed to find the zeros (and their multiplicity) of the expression

\begin{equation*}
    P(z)=\dfrac{z}{1!}+\dfrac{z^2}{2!}+\dfrac{z^3}{3!}+\cdots = e^z-1.
\end{equation*}
The general zero of this expression is $z=\log(1)$. But in his paper {\it ``De la controverse entre Mrs. Leibnitz et Bernoulli sur les logarithmes des nombres négatifs et imaginaires"} \cite{E168} (E168:``On the controverse of Leibniz and Bernoulli on the logarithms of negative and imaginary numbers") Euler had demonstrated that the logarithm of a number is a multivalued expression and hence concluded that there are infinitely many zeros, namely:

\begin{equation*}
    %\label{eq: Solutions Log}
    z= 0, \pm 2\pi i, \pm 4 \pi i, \pm 6 \pi i, \pm 8 \pi i, \cdots.
\end{equation*}
 Furthermore, all those zeros are simple, since:

\begin{equation*}
    \lim_{z \rightarrow 2k \pi i} \dfrac{e^z-1}{z- 2 k \pi i}=  \lim_{z \rightarrow 2k \pi i} \dfrac{e^z}{1} = e^{2  k \pi i}=1.
\end{equation*}
where L'Hospital's rule was used in the first step. Therefore, Euler  used the general solution formula (\ref{eq: Inhomogeneous Finite Solution}). This gave him:

\begin{equation}
    \label{eq: General Solution Simple Difference equation}
    f(x) = \int g(x)dx + e^{2 \pi i x}\int g(x) e^{-2 \pi i x}dx+  e^{-2 \pi i x}\int g(x) e^{+2 \pi i x}dx
\end{equation}
\begin{equation*}
    +  e^{4 \pi i x}\int g(x) e^{-4 \pi i x}dx + e^{-4 \pi i x}\int g(x) e^{+4 \pi i x}dx + \cdots
\end{equation*}
In \cite{E189} § 55, Euler expressed the solutions using sines and cosines instead of the exponentials that we used here. 

Thus, we arrived at Euler's general solution of the simple difference equation (\ref{eq: Simple Difference}). Unfortunately, as we will see below in section (\ref{subsec: Discussion of the Result}), there is a mistake in Euler's solution (\ref{eq: General Solution Simple Difference equation}).

\section{Application to the Factorial}
\label{sec: Application to the Factorial}

In \cite{E189} \S 56 - \S 60, Euler applied his general formula (\ref{eq: General Solution Simple Difference equation}) to the factorial\footnote{More precisely, Euler actually considered the difference equation satisfied by the $\Gamma$-function.}, i.e., the function $y(x)$ satisfying:

\begin{equation}
    \label{eq: Factorial}
    y(x+1)=x y(x).
\end{equation}
This equation can be transformed into a simple difference equation by taking logarithms. We have:

\begin{equation*}
    \log y(x+1) - \log y(x) = \log(x).
\end{equation*}

\subsection{Application of the General Formula}
\label{subsec: Application of the general Formula}

Applying (\ref{eq: General Solution Simple Difference equation}) with $f(x)= \log y(x)$ and $g(x)= \log(x)$ we get:

\begin{equation}
\label{eq: General Solution Factorial}
    f(x) = x\log x-x +C + e^{2 \pi i x}\int \log(x) e^{-2 \pi i x}dx+  e^{-2 \pi i x}\int\log(x) e^{+2 \pi i x}dx
\end{equation}
\begin{equation*}
    +  e^{4 \pi i x}\int \log(x) e^{-4 \pi i x}dx + e^{-4 \pi i x}\int \log(x) e^{+4 \pi i x}dx + \cdots
\end{equation*}
where $\int \log(x)dx$ was already evaluated and $C$ is a constant of integration{\footnote{This is  the solution Euler gave in \cite{E189} \S 59. But he represented his solution using sines and cosines.}.

\subsection{Derivation of Stirling's Formula}
\label{subsec: Derivation of Stirling's Formula}

\S 59 -\S 60 of \cite{E189} contain the derivation of Stirling's formula (\ref{eq: Stirling}) from (\ref{eq: General Solution Factorial}). Euler first evaluated the general expression:

\begin{equation*}
    e^{2 k \pi i x} \int e^{-2k \pi ix}\log (x)dx.
\end{equation*}
He did so by integrating by parts infinitely many times with $e^{-2 k \pi ix}$ as function to be integrated. In modern and compact notation the result is\footnote{Since Euler used $\sin(2k \pi x)$ and $\cos(2 \pi x)$ instead of $e^{-2 k \pi i x}$, his result differs from the one we will find. But the derivation is the same in both cases, of course.}:

\begin{equation*}
    e^{2 k \pi i x} \int e^{-2k \pi ix}\log (x)dx= -\dfrac{\log(x)}{2 k \pi i}+ \sum_{n=1}^{\infty} \dfrac{(-1)^n  (n-1)!}{(2k\pi i)^{n+1} x^n}+ C_k e^{2 k \pi i x}.
\end{equation*}
$C_k$ is a constant of integration. Proceeding in the same way for all other integrals, we have the formal identity:

\begin{equation*}
    \log y(x)= x\log x- x + C +\sum_{k\in \mathbb{Z} \setminus \lbrace 0 \rbrace} \left(C_k e^{2k \pi i x}-\dfrac{\log(x)}{2 k \pi i}+ \sum_{n=1}^{\infty} \dfrac{(-1)^n  (n-1)!}{(2k\pi i)^{n+1}x^n} \right)
\end{equation*}
Let us simplify the sum. First, we note that

\begin{equation*}
    C+ \sum_{k\in \mathbb{Z} \setminus \lbrace 0 \rbrace} C_k e^{2k \pi i x}  =: h(x)
\end{equation*}
is a general periodic function, i.e., it satisfies $h(x+1)=h(x)$ for all $x$. Next, 

\begin{equation*}
    \sum_{k\in \mathbb{Z} \setminus \lbrace 0 \rbrace} \dfrac{\log (x)}{2 k \pi i} =0,
\end{equation*}
since the terms cancel each other. Therefore, we just need to evaluate the double sum. By a formal calculation we have:

\begin{equation}
\label{eq: Double Sum}
\sum_{k\in \mathbb{Z} \setminus \lbrace 0 \rbrace}  \sum_{n=1}^{\infty} \dfrac{(-1)^n  (n-1)!}{(2k\pi i)^{n+1}x^n}   =  \sum_{n=0}^{\infty} \sum_{k=1}^{\infty} \dfrac{2}{k^{2n+2}} \cdot \dfrac{ (-1)^{n} (2n)!}{(2\pi)^{2n+2}\cdot x^{2n+1}}.
\end{equation}
The sum over $k$ had been evaluated by Euler. The general formula can found, e.g., in \cite{E130} and in modern notation reads:

\begin{equation}
    \sum_{k=1}^{\infty} \dfrac{1}{k^{2n}} = \dfrac{(-1)^{n-1} (2\pi)^{2n}B_{2n}}{2(2n)!},
\end{equation}
where $B_n$ is the $n$-th Bernoulli number. Inserting this into (\ref{eq: Double Sum}), we find:

\begin{equation*}
    \sum_{k\in \mathbb{Z} \setminus \lbrace 0 \rbrace}  \sum_{n=1}^{\infty} \dfrac{(-1)^n  (n-1)!}{(2k\pi i)^{n+1}x^n}= \sum_{n=0}^{\infty} 2 \cdot \dfrac{(-1)^{n} (2\pi)^{2n+2}B_{2n+2}}{2(2n+2)!}\cdot \dfrac{ (-1)^{n} (2n)!}{(2\pi)^{2n+2}\cdot x^{2n+1}}.
\end{equation*}
Many terms cancel such that:

\begin{equation*}
      \sum_{k\in \mathbb{Z} \setminus \lbrace 0 \rbrace}  \sum_{n=1}^{\infty} \dfrac{(-1)^n  (n-1)!}{(2k\pi i)^{n+1}x^n}= \sum_{n=1}^{\infty} \dfrac{B_{2n}}{(2n-1)2n x^{2n-1}}.
\end{equation*}
Therefore, inserting everything we found into (\ref{eq: General Solution Factorial}) we get:

\begin{equation}
\label{eq: Euler Factorial}
    \log y(x) = x\log x -x +h(x)+ \sum_{n=1}^{\infty} \dfrac{B_{2n}}{(2n-1)2n x^{2n-1}},
\end{equation}
where $h(x)$ satisfies $h(x+1)=h(x)$. This equation is to be understood as an asymptotic series  of course and is the formula Euler arrived at in \cite{E189} \S 60, Euler just substituted the explicit numbers for the Bernoulli numbers. 
Comparing (\ref{eq: Euler Factorial}) to (\ref{eq: Stirling}), the term $\log(\sqrt{2\pi})$ is still missing. In \cite{E189} Euler argued that it follows from considering a special case, e.g., $x=1$\footnote{More precisely, Euler argued that $h(x)$ is to be considered as constant in this case and the value of this constant is equal to the sum $1-\sum_{n=1}^{\infty} \frac{B_{2n}}{(2n-1)2n}$ which Euler claims to be $\frac{1}{2}\log (2\pi)$ without a proof in this paper, although the series does not converge due to the rapid growth of the Bernoulli numbers. But Euler knew that one can ascribe the beforementioned value to the sum, since it corresponds to the constant $\sqrt{2\pi}$ in Stirling's formula (\ref{eq: Stirling}).} and the initial condition $y(1)=1$ to (\ref{eq: Factorial}) such that one arrives at the final formula:

\begin{equation}
\label{eq: Euler Final log}
     \log y(x) = x\log x -x +\log (\sqrt{2 \pi})+ \sum_{n=1}^{\infty} \dfrac{B_{2n}}{(2n-1)2n x^{2n-1}},
\end{equation}
if $x$ is infinitely large. In \cite{E189} \S 60, Euler stated the formula as follows:

\begin{equation}
    \label{eq: Euler Factorial Final}
    y(x) = \dfrac{x^x}{e^x}\left(1+\dfrac{1}{12x}+\dfrac{1}{288x^2}-\dfrac{139}{51840x^3}+\cdots\right)\sqrt{2\pi},  
\end{equation}
which follows by inserting the explicit values for the Bernoulli numbers in (\ref{eq: Euler Final log}), taking the exponential and expanding the exponential of the sum.

\subsection{Discussion of the Result}
\label{subsec: Discussion of the Result}

As it was remarked by G. Faber in a footnote in the Opera Omnia version of \cite{E189}, equation (\ref{eq: Euler Final log}) and hence (\ref{eq: Euler Factorial Final}) is incorrect. The correct formula reads:

\begin{equation}
\label{eq: Log Stirling correct}
     \log y(x) = x\log x -x +\log (\sqrt{\dfrac{2 \pi}{x}})+ \sum_{n=1}^{\infty} \dfrac{B_{2n}}{(2n-1)2n x^{2n-1}},
\end{equation}
i.e., Euler's formula is off by the term $\log(\sqrt{x})$. Furthermore, the term is not missing due to a calculational error, but due to a conceptional one. More precisely, Euler's idea to construct the solution from the zeros of (\ref{eq:  Infinite Polynomial}) does not work in general.

We can see how the missing term enters by a formal argument\footnote{There are also rigorous arguments involving the theory of the Fourier transform. But this would carry us too far away from our actual objective.}. We are still interested in (\ref{eq: Simple Difference}). Writing $D$ for $\frac{d}{dx}$, this equation can also be represented as:

\begin{equation*}
    \left(e^D-1\right)f(x) =g(x).
\end{equation*}
Thus, formally the solution is given as:

\begin{equation*}
    f(x) =  \left(e^D-1\right)^{-1} g(x),
\end{equation*}
such that we have to find out how to express $\left(e^D-1\right)^{-1}$. We only know how to calculate $D^{n}f(x)$ for $n \in \mathbb{Z}$. Thus, the idea is to expand $\left(e^D-1\right)^{-1}$ into a Laurent series in $D$ around $D=0$ apply it to $g(x)$. There are many possibility to perform this expansion, but for purposes we will only need the direct expansion. This expansion had also been given by Euler, e.g., in {\it ``De seriebus quibusdam considerationes"} \cite{E130} (E130:``Considerations on certain series") \S 27\footnote{Euler considered the function $\frac{z}{1-e^{-z}}$ and did not state the general formula for the coefficients, but explained their origin.}. The expansion reads:

\begin{equation}
    %\label{eq: Expansion direct}
     \left(e^D-1\right)^{-1} = \sum_{n=0}^{\infty} B_n\dfrac{D^{n-1}}{n!}=D^{-1}-\dfrac{1}{2}+\dfrac{D}{12}- \dfrac{D^3}{720}+\cdots,
\end{equation}
where $B_n$ are the Bernoulli numbers again. Interpreting $D^{-1}$ as an integration, we can write:

\begin{equation}
\label{eq: Expansion Appliciation}
     f(x) =  \left(e^D-1\right)^{-1} g(x) = \int g(x)dx - \dfrac{1}{2}g(x) + \dfrac{1}{12}\dfrac{d}{dx}g(x)- \cdots,  
\end{equation}
which is nothing but a modern representation of the Euler-Maclaurin summation formula. Thus, Euler's approach, i.e., constructing the solution from the zeros of $e^D-1$, misses the term $-\frac{1}{2}g(x)$. If we apply (\ref{eq: Expansion Appliciation}) to the factorial, i.e., take $g(x)=\log(x)$
we arrive at (\ref{eq: Log Stirling correct}).

%%%%%%%%%%%%%%%%%%%%%%%%%%%%%%%%%%%%%%%%%%%%%%%%%%%%%%%%%%%%%%%%%%%%%%%%%%%%%%%%
\section{Conclusion}

In this note we briefly mentioned Euler's theory how to solve inhomogeneous ordinary differential equations of infinite order with constant coefficients and Euler's application of his theory to the derivation of Stirling's formula (\ref{eq: Stirling}). We pointed out the conceptual error in Euler's approach and provided an explanation how to correct it (section \ref{subsec: Discussion of the Result}). Nevertheless, there are many intriguing ideas in \cite{E189}, aside from Euler's derivation of Stirling's formula on which  we focused, such that we intend to cover more content from the before-mentioned paper in the future.

%%%%%%%%%%%%%%%%%%%%%%%%%%%%%%%%%%%%%%%%%%%%%%%%%%%%%%%%%%%%%%%%%%%%%%%%%%%%%%%%


\begin{thebibliography}{}




\bibitem[1]{E130}
Euler, L. (1750). ``De seriebus quibusdam considerationes" (E130). {\it Commentarii academiae scientiarum Petropolitanae, Volume 12}  (1739): pp.  53--96. Reprinted in Opera Omnia: Series 1, Volume 14, pp. 407--462. Original text available online at \url{https://scholarlycommons.pacific.edu/euler/}.

\bibitem[2]{E168}
Euler, L. (1751). ``De la controverse entre Mrs. Leibnitz et Bernoulli sur les logarithmes des nombres négatifs et imaginaires" (E168). {\it Mémoires de l'académie des sciences de Berlin, Volume 5}  (1747): pp.  139--179. Reprinted in Opera Omnia: Series 1, Volume 17, pp. 195--232. Original text available online at \url{https://scholarlycommons.pacific.edu/euler/}.

%\bibitem[3]{E170}
%Euler, L. (1751). ``Recherches sur les racines imaginaires des équations" (E170). {\it Mémoires de l'académie des sciences de Berlin, Volume 5}  (1746): pp.  222--288. Reprinted in Opera Omnia: Series 1, Volume 6, pp. 78--150. Original text available online at \url{https://scholarlycommons.pacific.edu/euler/}.

\bibitem[3]{E188}
Euler, L. (1753). ``Methodus aequationes differentiales altiorum graduum integrandi ulterius promota" (E188). {\it Novi Commentarii academiae scientiarum Petropolitanae, Volume 3}  (1750): pp.  pp. 3--35. Reprinted in Opera Omnia: Series 1, Volume 22, pp.181--213. Original text available online at \url{https://scholarlycommons.pacific.edu/euler/}.

\bibitem[4]{E189}
Euler, L. (1753). ``De serierum determinatione seu nova methodus inveniendi terminos generales serierum" (E189). {\it Novi Commentarii academiae scientiarum Petropolitanae, Volume 3} (1749): pp. 36--85. Reprinted in Opera Omnia: Series 1, Vol. 14, pp. 463--515. Original text available online at \url{https://scholarlycommons.pacific.edu/euler/}.



\end{thebibliography}
\end{document}